\documentclass{amsproc}
\usepackage{amsfonts}

\newtheorem{theorem}[subsection]{Theorem}

\newtheorem{proposition}[subsection]{Proposition}

\theoremstyle{definition}

\newtheorem{example}[subsection]{Example}

\theoremstyle{remark}
\newtheorem{remark}[subsection]{Remark}

\numberwithin{equation}{section}
\newcommand{\thref}[1]{Theorem \ref{#1}}
\newcommand{\prref}[1]{Proposition \ref{#1}}

\newcommand{\exref}[1]{Example \ref{#1}}

\newcommand\lbb[1]{\label{#1}}

\def\ti{\tilde}

\def\d{\partial}
\def\dd{{\mathrm{\, d}}}                           

\def\st{\; | \;}                               

\def\i{{\mathrm{i}}}     
\def\Cset{\mathbb{C}}       
\def\Rset{\mathbb{R}}       

\def\al{\alpha}                         

\def\Ga{\Gamma}
\def\de{\delta}
\def\De{\Delta}

\def\La{\Lambda}

\def\om{\omega}

\def\si{\sigma}

\def\th{\theta}

\def\A{{\mathcal{A}}}
\def\B{{\mathcal{B}}}
\def\K{{\mathcal{K}}}
\def\L{{\mathcal{L}}}

\DeclareMathOperator{\ad}{ad}

\begin{document}

\title[The Weyl algebra and bispectral operators]
{Automorphisms of the Weyl algebra and bispectral operators}

\author[B.~Bakalov]{Bojko Bakalov}
\address{Department of Mathematics, MIT, Cambridge, MA 02139, USA}
\email{bakalov@math.mit.edu}

\author[E.~Horozov]{Emil Horozov}
\address{Department of Mathematics and Informatics,
Sofia University, 5 J. Bourchier Blvd., Sofia 1126, Bulgaria}
\email{horozov@fmi.uni-sofia.bg}

\author[M.~Yakimov]{Milen Yakimov}
\address{Department of Mathematics, University of California,
Berkeley, CA 94720, USA}
\email{yakimov@math.berkeley.edu}

\date{July 14, 1997}


\maketitle
\section{Introduction}\lbb{s1}

The present note can be considered as an analytic counterpart of 
\cite{pl}.
In \cite{pl} we suggested several algebraic methods leading to new
(compared to \cite{DG, W1, Z, cmp, KR, GH1, GH2})
solutions of the bispectral problem in different versions.
Before sketching the main idea of \cite{pl}, we reformulate
the bispectral problem in a form which can be traced back to 
\cite{DG}
and which was used explicitly in \cite{pl, KR}.

Consider two associative algebras $\B$ and $\B'$ of differential 
operators.
We denote the variable in $\B$ by $x$ and the one in $\B'$ by $z$.
Let $\psi(x,z)$ be a function such that for any $P\in\B$
there exists $Q\in\B'$ satisfying
\begin{equation}\lbb{1.1}
P(x,\d_x)\psi(x,z) = Q(z,\d_z)\psi(x,z),
\end{equation}
and that for each $Q\in\B'$ there exists $P\in\B$ satisfying 
\eqref{1.1}.
We will assume that $P(x,\d_x)\psi(x,z) = 0$ implies $P=0$,
and similarly for $Q$. Then \eqref{1.1} defines an anti-isomorphism
\begin{equation}\lbb{b}
b\colon\B\to\B', \quad b(P)=Q,
\end{equation}
which depends on $\psi$.

We assume that both $\B$ and $\B'$ contain subalgebras of functions
$\K\subset\B$ and $\K'\subset\B'$.
Let $\A = b^{-1}(\K')$, $\A' = b(\K)$.
Then a solution to the {\em bispectral problem\/} of \cite{DG}
is a pair of operators of non-zero order $L\in\A$, $\La\in\A'$.
Indeed, if $b(L)=f\in\K'$ and  $b^{-1}(\La)=\th\in\K$,
then \eqref{1.1} implies
\begin{align}
\lbb{bisp1}
&L(x,\d_x)\psi(x,z) = f(z)\psi(x,z),
\\
\lbb{bisp2}
&\La(z,\d_z)\psi(x,z) = \th(x)\psi(x,z).
\end{align}

In \cite{pl} we pointed out that if $\si$ is an automorphism
of $\B$, then one can define a new anti-isomorphism
$b_1 = b\circ\si\colon\B\to\B'$. Then if the image
$b_1(\K)\subset\B'$ contains operators of non-zero order,
they will be bispectral. Of course, we need to define the
corresponding function $\psi_1(x,z)$ which determines $b_1$
via \eqref{1.1}.

The present note treats in detail the problem of defining $\psi_1$
in the case when both $\B$ and $\B'$ are two copies of the Weyl
algebra of regular differential operators in the complex plane 
$\Cset$.
This algebra is by now the most tractable example, as all its
automorphisms are known \cite{D}. They are all generated by the ones
of the form
\begin{equation*}
 e^{\ad p(x)} \;\;\text{or}\;\;  e^{\ad q(\d_x)}
\end{equation*}
where $p$ and $q$ are arbitrary polynomials (see \cite{D}).

We start with the obvious anti-isomorphism $b_0$ given by
\begin{equation}\lbb{1.2}
b_0(x)=\d_z, \quad b_0(\d_x)=z.
\end{equation}
The corresponding eigenfunction $\psi_0(x,z)$ satisfying
$P\psi_0 = b_0(P)\psi_0$, $P\in\B$, is $e^{xz}$.
Then any automorphism
\begin{equation}\lbb{2.1}
\si =  e^{\ad p_1(x)}  e^{\ad q_1(\d_x)}
\dotsm  e^{\ad p_m(x)}  e^{\ad q_m(\d_x)}
\end{equation}
of the Weyl algebra $\B$ gives an anti-isomorphism 
\begin{equation}\lbb{b2}
b=b_0\circ\si\colon\B\to\B'.
\end{equation}

Our objective is to make sense of the expressions of the type
\begin{equation}\lbb{1.3}
\psi(x,z) = e^{-q_m(\d_x)} e^{-p_m(x)} \dotsm e^{-q_1(\d_x)}
e^{-p_1(x)} e^{xz} \end{equation}
which will satisfy (\ref{1.1}, \ref{b}) for $b$. 
Indeed, we can write $\psi(x,z) = g e^{xz}$ with 
$g = e^{-q_m(\d_x)}\dotsm e^{-p_1(x)}$.
Then for $Q=Q(z,\d_z)\in\B'$, we have 
\begin{equation*}
b^{-1}(Q)\, \psi = \si^{-1}(b_0^{-1}Q) \,\psi 
= g\, (b_0^{-1}Q)\, g^{-1}\, g\, e^{xz}
= g\, (b_0^{-1}Q)\, e^{xz}
= g\, Q\, e^{xz} = Q\,\psi.
\end{equation*}
($g$ and $Q$ commute, since the former acts on $x$, while the latter on $z$.)

Our approach, as briefly 
indicated in \cite{pl}, is to apply repeatedly Laplace transformations
to the expression \eqref{1.3}. We will show that
for generic polynomials $p_j$, $q_j$, \eqref{1.3}
makes sense.
Moreover, the corresponding algebras $\A$ and $\A'$ contain 
differential
operators of non-zero order, for example, $L=b^{-1}(z)$ and 
$\La=b(x)$.
They satisfy
\begin{align}
\lbb{bisp3}
&L(x,\d_x)\psi(x,z) = z\psi(x,z),
\\
\lbb{bisp4}
&\La(z,\d_z)\psi(x,z) = x\psi(x,z),
\end{align}
hence they solve the bispectral problem.
It is easy to see that in this case
$\K=\Cset[x]$, $\K'=\Cset[z]$.

The main result of this note is \thref{t3.1} by which we
explicitly define $\psi(x,z)$ as a multiple Laplace integral.

\subsection*{Acknowledgements}
This paper is written after an
informal talk, one of us (E.H.) made at the Workshop on the Bispectral
Problem, March 1997, CRM, Montr\'eal, and on some of the material
from \cite{pl}. In many respects it is an answer to some
of the questions asked by the friendly audience.
We are particularly indebted to F.~A.~Gr\"unbaum and
G.~Wilson (to mention some of the names) who asked for this
talk. Also F.~A.~Gr\"unbaum persuaded us to write down the details. 
He kept interested in this work in the process of writing it
and himself made important remarks. 
We are grateful to him for all this. 
Last, but not least, we thank CRM and
Universit\'e de Montr\'eal, and in particular the organizers
J.~Harnad, A.~Kasman, and P.~Winternitz for the kind
hospitality and for inviting us to present our results.
E.H.\ acknowledges the partial support of Grant No.\ 523
of Bulgarian Ministry of Education.
B.B.\ acknowledges the hospitality of Universit\`a di Roma 
``La Sapienza'' during June 1997, when this work was completed.
\section{Special cases}\lbb{s2}
\subsection{}
In this section we will consider several
special cases of the construction outlined at the end of the previous
section.
The purpose of the analysis below is to
exclude from further study the sets of polynomials $p_j$, $q_j$
which do not produce new bispectral operators.

Let us first clarify which bispectral operators will be considered
non-trivially different. (We closely follow \cite{DG, W1}.)
First, by a change of the variable $x$, we can always make the 
leading
coefficient of $L$ constant. Then multiplying the eigenfunction 
$\psi$
by a function depending only on $x$, we can make the second
coefficient of $L$ zero. This fixes the operator $L$ up to an
affine change $x\mapsto ax+b$. In the same manner, we fix $\La$.
As proved in \cite{DG, W1}, in this case both operators $L$ and $\La$
have rational coefficients.

Thus, we will consider only bispectral operators with constant
leading coefficients and vanishing second ones;
any such operator
defined up to affine changes of variables.
Nevertheless, in some of the intermediate
steps we may use operators which differ by the above mentioned
trivial transformations.

\subsection{}
First, it is easy to see that all polynomials $p_j$, $q_j$
should have positive degrees. Otherwise, using the above 
transformations,
one can reduce $\psi$ from \eqref{1.3} to a function of this kind
but with smaller $m$.

Indeed, this is obvious for $q_1,p_2,\dots,q_{m-1},p_m,q_m$.
Suppose that $\deg p_1=0$. Then we can exchange the roles of $x$ and $z$,
considering 
\begin{equation*}
\ti\psi(x,z) := \psi(z,x) = 
e^{-p_1(\d_x)} e^{-q_1(x)} \dotsm e^{-p_m(\d_x)}
e^{-q_m(x)} e^{xz}.
\end{equation*}

\subsection{}
If $\deg p_1 = 1$, this amounts to an affine change of the variable
$x$ (cf.\ \eqref{1.3}). Exchanging the roles of $x$ and $z$,
as above, we can exclude the case when $\deg q_1 = 1$.
The same argument can be used to eliminate any polynomial
$p_j$ or $q_j$ of degree $1$. (This amounts to  
affine changes of the variables $u_j$, $v_j$ from 
\eqref{3.4} below.)

Therefore, from now on we will consider only polynomials
$p_j$, $q_j$ of degree at least $2$.

\begin{example}\lbb{e2.2}
Let $\si = e^{\ad p(x)} \,  e^{\ad q(\d_x)}$.
Then, as computed in \cite{pl}, we have
\begin{align*}
&b(x) = b_0\si(x) = \d_z + q'(z-p'(\d_z)),
\\
&b(\d_x) = b_0\si(\d_x) = z-p'(\d_z).
\end{align*}
We will prove in the next section that the function
$\psi(x,z) = e^{-q(\d_x)} e^{-p(x)} e^{xz}$ exists.
Thus we find the bispectral operators
\begin{align}
\lbb{2.3}
&L(x,\d_x) = b^{-1}(z) = \d_x + p'(x-q'(\d_x)), 
\\
\lbb{2.4}
&\La(z,\d_z) = b(x) = \d_z + q'(z-p'(\d_z)), 
\end{align}
satisfying (\ref{bisp3}, \ref{bisp4}).

Suppose that $\deg p=2$, i.e.\ that $p(x)=a x^2 + b x$.
Then
\begin{align*}
&L(x,\d_x) = \d_x - 2 a q'(\d_x) + 2 a x + b,
\\
&\La(z,\d_z) = \d_z + q'(z - 2 a \d_z - b).
\end{align*}
By affine changes of the variables $x$ and $z$
and a multiplication of $\psi(x,z)$ by $e^{p(x)}$,
it is easy to see that
$L$ and $\La$ are the (generalized) Airy operators (see \cite{cmp, 
KR})
provided that $\deg q \ge3$. The same result holds if $\deg p\ge3$, 
$\deg q=2$.

If both $\deg p\ge3$, $\deg q \ge3$, we obtain a new bispectral pair
of operators $L$, $\La$, as explained in \cite{pl}.

The case when $\deg p=2$, $\deg q=2$ is trivial and is treated in
the next subsection.
\end{example}

\subsection{}
Now suppose that all polynomials $p_j$, $q_j$ have degree $2$.
Then one easily gets (by induction) that
\begin{align*}
&b(x) = b_0\si(x) = a_{11} z + a_{12} \d_z,
\\
&b(\d_x) = b_0\si(\d_x) = a_{21} z + a_{22} \d_z
\end{align*}
for some non-degenerate matrix $a_{ij}$.

If $a_{12}=0$, then the image of $\Cset[x]$ under $b$ is 
$\Cset[z]$,
i.e.\ $b$ does not produce bispectral operators.
In the case when $a_{12}\ne0$, $b$
gives bispectral operators of rank $1$. (For the notion of rank,
cf.\ \cite{DG, W1, cmp}).

\section{The existence of joint eigenfunctions}\lbb{s3}
\subsection{}
In this section we will demonstrate the existence of a common 
eigenfunction
$\psi(x,z)$ of (\ref{bisp3}, \ref{bisp4}). We will construct the 
function
$\psi(x,z)$ by a repeated application
of Laplace transformation on \eqref{1.3} and induction on
$m$. According to Sect.~\ref{s2},
it suffices to consider polynomials $p_j$, $q_j$ of degree at least 
$2$, with one of them of degree greater than $2$.
\subsection{}

Our starting point is the function $\psi_0(x,z) = e^{xz}$.
%
%
The crucial step is to define
$\psi_1(x,z) = e^{-q_1(\d_x)} e^{-p_1(x)} e^{xz}$.
Let us first explain the idea of the construction
as it  will be applied later for any $m$.
For simplicity, we will drop the subscripts $1$.

Supposing that $\psi(x,z)$ exists, we can perform a
Laplace transformation on it with respect to $x$:
\begin{multline*}
(\L\psi)(v,z) = \int_{\Ga'} \psi(u,z) e^{-uv} \dd u
\\
= \int_{\Ga'} \left( e^{-q(\d_u)} e^{-p(u)} e^{uz} \right) e^{-uv} \dd u
= e^{-q(v)} \int_{\Ga'} e^{-p(u) + uz - uv} \dd u.
\end{multline*}
We choose the contour $\Ga'$ so that the last integral makes
sense. This can be done, for example, as follows.

Let $n=\deg p\ge2$ and fix an $n$-th root $\al$ of the leading
coefficient of $p$. Let $\om_1$, $\om_2$ be two different $n$-th roots
of unity. Define the contour $\Ga'$ to be the
union of two rays
\begin{equation}
\lbb{ga1'}
\Ga' = \{ \al^{-1} \om_1 t \st t\in(\infty,0] \}
\cup \{ \al^{-1} \om_2 t \st t\in[0,\infty) \}.
\end{equation}

Then the Laplace transformation
\begin{multline*}
(\L\psi)(v,z) = e^{-q(v)} \al^{-1} \om_1 \int_0^\infty
e^{-t^n+\dotsm} e^{(z-v) \al^{-1} \om_1 t} \dd t
\\
+ e^{-q(v)} \al^{-1} \om_2 \int_0^\infty
e^{-t^n+\dotsm} e^{(z-v) \al^{-1} \om_2 t} \dd t
\end{multline*}
is well defined,
i.e.\ the integrals are convergent, because the leading terms in the
exponentials are negative.

In the same manner, we can find the inverse Laplace transformation
of $(\L\psi)(v,z)$, i.e.\  to define $\psi(x,z)$.
Below, for simplicity, we drop the inessential factors of the type
$2\pi$. We define a contour $\Ga''$ using the above recipe for the
polynomial $q(v)$. Then we define
\begin{equation}\lbb{3.3}
\psi(x,z) = \int\limits_{\Ga''} (\L\psi)(v,z) e^{xv} \dd v
= \int\limits_{\Ga'}\int\limits_{\Ga''}
e^{-p(u) - q(v) - uv + xv + uz} \dd u \dd v.
\end{equation}
\begin{proposition}\lbb{psi1}
Let $p$ and  $q$ be polynomials of degree at least $2$. Define the
contours $\Ga'$ and $\Ga''$ by \eqref{ga1'} for $p$ and $q$, respectively.
Then, provided that one of the polynomials $p$, $q$ has degree at least $3$,
the function $\psi(x,z)$ from \eqref{3.3} is well defined and satisfies 
\textup{(\ref{bisp3}, \ref{bisp4})}  for the operators $L$ and $\La$ from
\textup{(\ref{2.3}, \ref{2.4})}.
\end{proposition}
\begin{proof}
(i) Let us first prove (\ref{bisp3}, \ref{bisp4}), assuming that
all integrals are absolutely convergent. For brevity, denote the 
integrand of \eqref{3.3} by $exp$. By differentiating under the integral,
we obtain
\begin{align*}
&\d_x\psi(x,z) = \iint v \, exp\, \dd u \dd v,
\\
&\d_z\psi(x,z) = \iint u \, exp\, \dd u \dd v.
\\
\intertext{Integration by parts gives}
&0 = \iint \d_u(exp) \dd u \dd v 
= \iint\bigl(-p'(u)-v+z\bigr)\, exp\,\dd u \dd v,
\\
&0 = \iint \d_v(exp) \dd u \dd v 
= \iint \bigl(-q'(v)-u+x\bigr)\, exp\,\dd u \dd v.
\end{align*}
Then
\begin{multline*}
z\psi(x,z) = \iint \bigl( p'(u)+v \bigr)\, exp\, \dd u \dd v 
= \iint \Bigl( p'\bigl(x-q'(v)\bigr)+v \Bigr)\, exp\, \dd u \dd v
\\
= \iint \Bigl( p'\bigl(x-q'(\d_x)\bigr)+\d_x \Bigr)\, exp\, \dd u \dd v
= \Bigl( p'\bigl(x-q'(\d_x)\bigr)+\d_x \Bigr) \psi(x,z).
\end{multline*}
The equation \eqref{bisp4} is proved similarly.

(ii) We have already explained that the integrals are convergent when
$\deg p\ge3$, $\deg q\ge3$. If one of the polynomials, say $p$, 
has degree $2$,
then by an affine change $u\mapsto u +  c v$ $(c\in\Cset)$
one can bring $-p(u) - uv$ to the form $\al u^2 + \mu v^2$.
Then on the integration domain the leading terms in the exponential
in \eqref{3.3} are negative. This proves the absolute convergence 
of the integrals.
\end{proof}
\begin{remark}
If both polynomials $p$, $q$
were of degree $2$, their quadratic part might have rank $1$,
i.e.\ it might be a perfect square. Then  $\psi(x,z)$ would
be divergent; in fact, it would be of a $\de$-function type.
\end{remark}

\subsection{}

The general case
is treated in the same manner. We are looking for
an eigenfunction $\psi(x,z)$ of the form
\begin{multline}\label{3.4}
\psi_m(x,z) = \int\limits_{\Ga'_1} \int\limits_{\Ga''_1} 
\dotsi \int\limits_{\Ga'_m} \int\limits_{\Ga''_m}
\exp\left(
\sum_{s=0}^m (u_{s+1}-u_s) v_s 
- p_s(u_s) - q_s(v_s) \right) 
\\
\dd u_1 \dd v_1 \dotsm \dd u_m \dd v_m ,
\end{multline}
where $u_{m+1} := x$, $q_0=p_0=u_0 :=0$, $v_0 := z$.

When all polynomials $p_1, q_1, \dots,
p_m, q_m$  have degrees greater than $3$,
the above explained choice of the contours
$\Ga'_1, \Ga''_1, \dots, \Ga'_m, \Ga''_m$ will automatically  make the integral
convergent.
The difficulties start when some of the polynomials have degrees equal to $2$.
The same argument as above works when 
{\em{in the sequence of polynomials $p_1, q_1, \dots,
p_m, q_m$  there are no two consecutive polynomials of
degree $2$}}.~{\footnote{
One can prove the convergence under weaker assumptions; 
this is more complicated and will be done elsewhere.}}
\begin{theorem}\lbb{t3.1}
Under the above assumptions, the integral \eqref{3.4} 
is convergent absolutely and uniformly in any compact. 
It satisfies \textup{(\ref{1.1}, \ref{b})} with $b$ given
by \textup{(\ref{2.1}, \ref{b2})}.
In particular, it solves the
bispectral problem \textup{(\ref{bisp3}, \ref{bisp4})}
for $L=b^{-1}(z)$, $\La=b(x)$.
\end{theorem}
\begin{proof}
is very similar to that of \prref{psi1} and is by induction on $m$.
Let us first introduce a shorthand notation: 
\begin{equation*}
L_m := b_m^{-1}(z), \quad D_m := b_m^{-1}(\d_z), \quad
\La_m := b_m(x), \quad \De_m := b_m(\d_x).
\end{equation*} 
We rewrite \eqref{3.4} as
\begin{equation}\lbb{3.5}
\psi_m(x,z) = \int\limits_{\Ga'_m} \int\limits_{\Ga''_m}
e^{ (x-u_m) v_m - p_m(u_m) - q_m(v_m) } \,\psi_{m-1}(u_m,z)
\dd u_m \dd v_m,
\end{equation}
or shortly
\begin{equation*}
\psi_m(x,z) = \iint e_m \psi_{m-1}.
\end{equation*}
By differentiating under the integral, we have
\begin{align}
\lbb{p3}
&\d_x \psi_m(x,z) = \iint v_m e_m \psi_{m-1}.
\\
\intertext{By the inductive assumption}
\lbb{p4}
&\La_{m-1}(z,\d_z) \,\psi_m(x,z) = \iint u_m e_m \psi_{m-1}.
\\
\intertext{Integration by parts in \eqref{3.5} gives}
\lbb{p1}
&0 = \iint e_m \,\bigl( -p'_m(u_m) - v_m + \De_{m-1} \bigr)\,\psi_{m-1},
\\
\lbb{p2}
&0 = \iint e_m \,\bigl( -q'_m(v_m) - u_m + x \bigr)\,\psi_{m-1}.
\end{align}
(In \eqref{p1} we again used the inductive assumption.)
For later use, we note that \eqref{p2} remains valid if we replace
$\psi_{m-1}$ with any other function independent of $v_m$.

Now we have
\begin{multline}\lbb{dem}
\d_x \psi_m(x,z) = \iint v_m e_m \psi_{m-1} 
\\
= \iint e_m \,\bigl( -p'_m(u_m) + \De_{m-1} \bigr)\,\psi_{m-1}
= \Bigl( \De_{m-1}(z,\d_z) - p'_m\bigl(\La_{m-1}(z,\d_z)\bigr) 
\Bigr) \,\psi_m(x,z),
\end{multline}
and similarly,
\begin{multline}\lbb{lam}
x \psi_m(x,z) = \iint e_m \,\bigl( q'_m(v_m) + u_m \bigr)\,\psi_{m-1}
\\
= \bigl( q'_m(\d_x) + \La_{m-1}(z,\d_z) \bigr)\,\psi_m(x,z)
\\
= \Bigl( q'_m \bigl(\De_{m-1}(z,\d_z) - p'_m\bigl(\La_{m-1}(z,\d_z)\bigr) 
\bigr) + \La_{m-1}(z,\d_z) \Bigr)\,\psi_m(x,z).
\end{multline}
To show that the operator in the last equality of \eqref{dem} is
exactly $\De_m$, we note that 
\begin{equation*}
\De_m = b_m(\d_x) = b_{m-1}(\si_m(\d_x)) 
= b_{m-1} \bigl( \d_x - p'_m(x) \bigr)
= \De_{m-1} - p'_m(\La_{m-1}).
\end{equation*}
Here $\si_m := e^{\ad p_m(x)} \,  e^{\ad q_m(\d_x)}$, cf.\ \exref{e2.2}.
In the same way, \eqref{lam} implies $x \psi_m = b_m(x) \psi_m$.

Using the inductive assumption, \eqref{p2} and integration by parts, we compute
\begin{multline*}
z \psi_m(x,z) = \iint e_m L_{m-1}(u_m, \d_{u_m})\,\psi_{m-1}
\\
= \iint e_m L_{m-1}\bigl( x-q'_m(v_m), \d_{u_m} \bigr)\,\psi_{m-1}
= \iint \Bigl( L_{m-1}\bigl( x-q'_m(v_m), -\d_{u_m} \bigr) e_m \Bigr)
\,\psi_{m-1}
\\
= \iint L_{m-1}\bigl( x-q'_m(v_m), p'_m(u_m) + v_m \bigr) \, e_m \,\psi_{m-1}
\\
= L_{m-1}\Bigl( x-q'_m(\d_x), p'_m\bigl(x-q'_m(\d_x)\bigr) + \d_x \Bigr) 
\,\psi_m(x,z).
\end{multline*}
The last operator is exactly 
\begin{equation*}
L_{m-1}\bigl( \si_m^{-1}(x), \si_m^{-1}(\d_x) \bigr)
= \si_m^{-1} \bigl(L_{m-1}(x,\d_x)\bigr)  
= \si_m^{-1} b_{m-1}^{-1}(z) = b_m^{-1}(z).
\end{equation*}

The proof that $\d_z \psi_m(x,z) = D_m(x,\d_x)\psi_m(x,z)$ is similar
and is left to the reader.
\end{proof}
\subsection{}
Here we will present some peculiar properties of the above solutions,
discovered by F.~A.~Gr\"unbaum  \cite{G},
which show that some ``obvious'' conjectures concerning the
bispectral problem, one is tempted to pose, are not true.

Let $m=1$, $p_1(x) = q_1(x) = x^3 / 3$
and define contours
\begin{equation*}
\Ga_k = \Rset_+ \cup 
\{ e^{2\pi \i k / 3} t \st t\in [0,\infty) \} ,
\qquad k=1,2
\end{equation*}
with arbitrary directions.

Then, by \prref{psi1}, the functions
\begin{equation*}
\psi_{kl}(x,z) = \int\limits_{\Ga_k}\int\limits_{\Ga_l}
e^{-u^3/3 - v^3/3 -uv + xv + uz} \dd u \dd v
\end{equation*}
satisfy the equations 
\begin{align*}
&x \psi_{kl}(x,z) =
\left(\d_z +\left(z-\d_z^2 \right)^2 \right) \psi_{kl}(x,z),
\\
&z \psi_{kl}(x,z) =
\left(\d_x +\left(x-\d_x^2 \right)^2 \right)  \psi_{kl}(x,z).
\end{align*}

Note that the operators on the right hand sides are
identical (with $x$ and $z$ switched). As it is well
known, in the Bessel case or Airy case 
all joint eigenfunctions
are symmetric with respect to $x$ and $z$  (see \cite{DG, cmp, KR}).
But as it was noticed by F.~A.~Gr\"unbaum
\cite{G}, here the situation is different.
For example
$\psi_{11}$ and $\psi_{22}$ are
obviously symmetric, but $\psi_{12}$ and $\psi_{21}$
are not. Only their sum is a symmetric function. In this way
we give another proof of Gr\"unbaum's result that the space
of symmetric eigenfunctions is only 3-dimensional.

\end{document}